\documentclass[leqno,draft]{article}

\newtheorem{theorem}{Theorem}
\newtheorem{lemma}[theorem]{Lemma}
\newtheorem{proposition}[theorem]{Proposition}
\newtheorem{definition}[theorem]{Definition}
\newtheorem{corollary}[theorem]{Corollary}

\newcommand{\begintheorem}{\addtocounter{equation}{1}\begin{theorem}}
\newcommand{\beginlemma}{\addtocounter{equation}{1}\begin{lemma}}
\newcommand{\beginproposition}{\addtocounter{equation}{1}\begin{proposition}}
\newcommand{\begindefinition}{\addtocounter{equation}{1}\begin{definition}}
\newcommand{\begincorollary}{\addtocounter{equation}{1}\begin{corollary}}

\begin{document}

\title{Calculus, fractals, and analysis on metric spaces}

\author{Stephen Semmes\thanks{This article has been prepared in connection
with the 2009 Barnett lecture at the University of Cincinnati, and the author
is grateful to the mathematics department there for the kind invitation.} \\
        Rice University}

\date{}

\maketitle

\centerline{\it Dedicated to the memory of Juha Heinonen}

\begin{abstract}
A consequence of the mean value theorem is that a differentiable
function $f$ on an open interval $I$ is constant if $f'(x) = 0$ for
each $x \in I$.  This does not work so well in the context of fractal
sets, even when they are connected, as in a classical example of
Whitney.  More precisely, the behavior of functions on a set depends
on the geometry of the set, and there are numerous possibilities.
Here we discuss some basic notions and examples related to these
themes.
\end{abstract}

\tableofcontents

\section{Whitney's example}
\label{whitney's example}
\setcounter{equation}{0}

        If $f$ is a differentiable real-valued function on the real
line, $I \subseteq {\bf R}$ is an interval, and
\begin{equation}
        f'(x) = 0
\end{equation}
for every $x \in I$, then $f$ is constant on $I$, by the mean value
theorem.  Of course, there are analogous results for connected open
sets or smooth curves and surfaces in higher dimensions.  However, a
well-known example due to Whitney \cite{hw2} shows that a
continuously-differentiable function $f$ on the plane may satisfy
\begin{equation}
        \nabla f(p) = 0
\end{equation}
for every point $p$ in a connected set $E$ and not be constant on $E$.  

        More precisely, $E$ is a continuous curve in ${\bf R}^2$ in
Whitney's example, which cannot be very smooth.  At the same time, if
$f$ were too smooth, then Whitney's example would also not work.  For
instance, if $f$ were twice-continuously differentiable, then $f(E)$
would have $1$-dimensional measure $0$ in the real line, by Sard's
theorem.  This would imply that $f(E)$ contains only one element,
since $f(E)$ is connected, by the connectedness of $E$ and continuity
of $f$.  This is the same as saying that $f$ is constant on $E$.

        Let us briefly review the proof of Sard's theorem in this
context.  If $p \in E$ and the second derivatives of $f$ are equal to
$0$ at $p$, then $f$ is very flat at $p$, by Taylor's theorem.  In
this case, it does not matter how complicated $E$ might be near $p$.
Otherwise, if one of the second derivatives of $f$ is not $0$ at $p$,
then the implicit function theorem implies that the set where one of
the first derivatives of $f$ is equal to $0$ is a
continuously-differentiable curve near $p$.  Thus one can use the
vanishing of the gradient directly on this part.

        One of Whitney's main tools was his celebrated extension
theorem \cite{hw1} for smooth functions.  This permits one to focus on
functions on $E$, which can be extended to ${\bf R}^2$ afterward.  For
that matter, we can consider functions on abstract metric spaces.

\section{Metric spaces}
\label{metric spaces}
\setcounter{equation}{0}

        Let $(M, d(x, y))$ be a metric space.  Thus $d(x, y)$ is a
nonnegative real-valued function on the set $M$, $d(x, y) = 0$
exactly when $x = y$,
\begin{equation}
        d(y, x) = d(x, y)
\end{equation}
for every $x, y \in M$, and
\begin{equation}
        d(x, z) \le d(x, y) + d(y, z)
\end{equation}
for every $x, y, z \in M$.  Of course, the latter is known as the
\emph{triangle inequality}, and ${\bf R}^n$ is a metric space with the
standard Euclidean metric.

        There is an obvious way to say what it means for a real-valued
function $f$ on $M$ to have gradient at a point $p \in M$ equal to
$0$, which is that
\begin{equation}
\label{lim_{x to p} frac{|f(x) - f(p)|}{d(x, p)} = 0}
        \lim_{x \to p} \frac{|f(x) - f(p)|}{d(x, p)} = 0.
\end{equation}
If $M$ is ${\bf R}^n$ with the standard Euclidean metric, then this is
equivalent to $f$ being differentiable at $p$, and with first
derivatives equal to $0$ at $p$.

        If (\ref{lim_{x to p} frac{|f(x) - f(p)|}{d(x, p)} = 0}) holds
at every point $p$ in a closed set $E \subseteq M$, then one may wish
to ask that the limit be uniform on compact subsets of $E$.  In particular,
this holds for continuously-differentiable functions on ${\bf R}^n$.
Conversely, if $E$ is a closed set in ${\bf R}^n$, and $f$ is a real-valued
function on $E$ that satisfies (\ref{lim_{x to p} frac{|f(x) - f(p)|}{d(x, p)} 
= 0}) uniformly on compact subsets of $E$, then Whitney's extension theorem
implies that $f$ can be extended to a continuously-differentiable function
on ${\bf R}^n$ whose gradient is equal to $0$ on $E$.

        A real-valued function $f$ on a metric space $M$ is said to be
\emph{Lipchitz} if there is a nonnegative real number $C$ such that
\begin{equation}
\label{|f(x) - f(y)| le C d(x, y)}
        |f(x) - f(y)| \le C \, d(x, y)
\end{equation}
for every $x, y \in M$.  Actually, it suffices to check that
\begin{equation}
        f(x) \le f(y) + C \, d(x, y)
\end{equation}
for every $x, y \in M$, since this is the same as
\begin{equation}
        f(x) - f(y) \le C \, d(x, y),
\end{equation}
and it implies that
\begin{equation}
        f(y) - f(x) \le C \, d(x, y)
\end{equation}
by interchanging $x$ and $y$.  For example, $f_a(x) = d(a, x)$ is
Lipschitz with $C = 1$ for each $a \in M$, by the triangle inequality.

\section{Snowflake spaces}
\label{snowflake spaces}
\setcounter{equation}{0}

        If $a$, $b$ are nonnegative real numbers and $0 < t < 1$, then
\begin{equation}
\label{(a + b)^t le a^t + b^t}
        (a + b)^t \le a^t + b^t.
\end{equation}
For example, if $t = 1/2$, then this can be seen by squaring both sides.
Otherwise, one can observe that
\begin{equation}
        \max(a, b) \le (a^t + b^t)^{1/t},
\end{equation}
and hence
\begin{eqnarray}
        a + b & \le & (a^t + b^t) \, \max(a, b)^{1 - t}                  \\
              & \le & (a^t + b^t) \, (a^t + b^t)^{(1 - t)/ t}  \nonumber \\
               & = & (a^t + b^t)^{1/t},                         \nonumber
\end{eqnarray}
which is equivalent to (\ref{(a + b)^t le a^t + b^t}).  One can also
use calculus to show (\ref{(a + b)^t le a^t + b^t}).

        If $(M, d(x, y))$ is any metric space, then it follows that
$d(x, y)^t$ is a metric on $M$ when $0 < t < 1$ too.  The main point
is that $d(x, y)^t$ satisfies the triangle inequality, because of
(\ref{(a + b)^t le a^t + b^t}).  Note that $d(x, y)^t$ determines the
same topology on $M$ as $d(x, y)$.

        Suppose that $f$ is a Lipschitz function on $M$ with respect
to $d(x, y)$.  It is easy to see that
\begin{equation}
        \lim_{x \to p} \frac{|f(x) - f(p)|}{d(x, p)^t} = 0
\end{equation}
uniformly on $M$ when $0 < t < 1$.  Hence there are plenty of nonconstant
functions on $M$ with ``gradient $0$'' with respect to $d(x, y)^t$.

        In particular, this can be applied to connected metric spaces.
For instance, there are plenty of nonconstant functions with
``gradient $0$'' on the unit interval $[0, 1]$ with respect to the
metric $|x - y|^t$ when $0 < t < 1$.

\section{Bilipschitz embeddings}
\label{bilipschitz embeddings}
\setcounter{equation}{0}

        Let $(M_1, d_1(x, y))$ and $(M_2, d_2(u, v))$ be metric
spaces.  A mapping $\phi$ from $M_1$ into $M_2$ is said to be
\emph{bilipschitz} if there is a real number $C \ge 1$ such that
\begin{equation}
\label{C^{-1} d_1(x, y) le d_2(phi(x), phi(y)) le C d_1(x, y)}
        C^{-1} \, d_1(x, y) \le d_2(\phi(x), \phi(y)) \le C \, d_1(x, y)
\end{equation}
for every $x, y \in M_1$.  Equivalently, $\phi$ is bilipschitz if it
is one-to-one, Lipschitz, and its inverse is Lipschitz as a mapping
from $\phi(M_1)$ onto $M_1$.

        Many standard snowflake curves in the plane are
bilipschitz-equivalent, at least locally, to an interval in the real
line with a snowflake metric $|x - y|^t$ for some $0 < t < 1$.  For
instance, the von Koch snowflake curve has this property with $t =
\log 3 / \log 4$.

        If $\phi : M_1 \to M_2$ is bilipschitz, then it is easy to see
that a real-valued function $f$ on $\phi(M_1)$ has ``gradient $0$'' if
and only if $f \circ \phi$ has ``gradient $0$'' on $M_1$.  Thus one
can get examples of nonconstant functions with ``gradient $0$'' on
snowflake curves, using the corresponding examples on the real line
with respect to $|x - y|^t$ when $0 < t < 1$.  As mentioned
previously, one can then use Whitney's extension theorem to get
continuously-differentiable functions on the plane with the desired
properties.  One can also look at snowflake curves in terms of
quasiconformal mappings on the plane.  See \cite{ph2} for more
sophisticated examples using snowflake metrics and embeddings, and
\cite{n1, n2} for some other results related to Whitney's example.

        Whitney's original example was based on Cantor sets instead of
snowflake curves, but the gist of the argument is quite similar.  Of
course, it is easy to find locally-constant functions on Cantor sets
that are not constant, because of disconnectedness.  It is more
interesting to consider mappings from Cantor sets onto intervals of
positive length, and which are as smooth as possible.  For instance,
there is a continuous function from the middle-thirds Cantor set in
the unit interval $[0, 1]$ onto $[0, 1]$.  Using a Cantor set in the
plane of Hausdorff dimension $1$, one can get a Lipschitz function
onto $[0, 1]$.  With a slightly larger Cantor set in the plane, one
can get a function that has ``gradient $0$'' and maps onto $[0, 1]$.
To get a connected set, one can pass a curve through the Cantor set in
an appropriate manner.

\section{Nice curves}
\label{nice curves}
\setcounter{equation}{0}

        Let $(M, d(x, y))$ be a metric space.  As usual, a continuous
path in $M$ is defined by a continuous mapping from a closed interval
$[a, b]$ in the real line into $M$.  If $q : [b, c] \to M$ is another
continuous path and $p(b) = q(b)$, then we can combine the two in the
obvious way to get a continuous path defined on $[a, c]$.

        Let us say that $p : [a, b] \to M$ is a \emph{nice curve} if
\begin{equation}
\label{d(p(x), p(y)) le |x - y|}
        d(p(x), p(y)) \le |x - y|
\end{equation}
for every $x, y \in [a, b]$, which is the same as saying that $p$ is
Lipschitz with constant $1$ as a mapping from $[a, b]$ into $M$.  If
$q : [b, c]$ is another nice curve with $p(b) = q(b)$, then the
combined path on $[a, c]$ is nice too.

        A continuous path $p : [a, b] \to M$ is said to have
\emph{finite length} if the sums
\begin{equation}
\label{sum_{i = 1}^n d(p(x_i), p(x_{i + 1}))}
        \sum_{i = 1}^n d(p(x_i), p(x_{i + 1}))
\end{equation}
are bounded, where
\begin{equation}
        a = x_0 < x_1 < x_2 \cdots < x_n = b
\end{equation}
is any partition of $[a, b]$.  In this case, the \emph{length} of the
path is defined to be the supremum of (\ref{sum_{i = 1}^n d(p(x_i),
p(x_{i + 1}))}) over all partitions of $[a, b]$.  If $p$ is nice,
then $p$ has finite length less than or equal to $b - a$.

        Conversely, if $p$ is a continuous path of finite length, then
there is a reparameterization of $p$ which is a nice curve defined on
an interval whose length is the length of $p$.  Basically, one can
reparameterize $p$ by arc length, as in vector calculus.  For
simplicity, we shall restrict our attention to nice curves here.

\section{Calculus}
\label{calculus}
\setcounter{equation}{0}

        Let $(M, d(x, y))$ be a metric space, and let $p : [a, b] \to
M$ be a nice curve in $M$.  Also let $f$ be a real-valued function on
$M$ which has ``gradient $0$'' in the sense of (\ref{lim_{x to p}
frac{|f(x) - f(p)|}{d(x, p)} = 0}) at $p(x)$ for every $x \in [a, b]$.
This implies that $f(p(x))$ has derivative $0$ in the usual sense of
calculus when $a < x < b$, with one-sided derivatives equal to $0$ at
the endpoints.  In particular, $f(p(x))$ is continuous on $[a, b]$,
and the mean value theorem implies that $f(p(x))$ is constant on $[a,
b]$.

        Alternatively, one might start with the hypotheses that $p :
[a, b] \to M$ be a continuous path of finite length, and that
(\ref{lim_{x to p} frac{|f(x) - f(p)|}{d(x, p)} = 0}) hold with
uniform convergence on $p([a, b])$.  In this case, one can show
somewhat more directly that $f(p(x))$ is constant on $[a, b]$,
basically by showing that it has length $0$ as a continuous path in
the real line.  It suffices to consider partitions of $[a, b]$ in
which the distances between consecutive elements is small, because the
relevant sums can only get larger as points are added to a partition,
by the triangle inequality.  The sums for $f(p(x))$ can then be
estimated in terms of the gradient $0$ condition for $f$ and the
finite length of $p$.

        Of course, there are many variants of arguments like these.
Well-known theorems in real analysis state that Lipschitz functions on
the real line are differentiable almost everywhere with respect to
Lebesgue measure, and that the fundamental theorem of calculus holds
for Lipschitz function using the Lebesgue integral of the derivative.
Thus a Lipschitz function on ${\bf R}$, or on an interval, is constant
when its derivative is equal to $0$ almost everywhere.  There are also
analogous results on ${\bf R}^n$.

        Monotone functions on the real line are differentiable almost
everywhere too.  This applies as well to functions of ``bounded
variation'', which amounts to finite length for paths in ${\bf R}$.
Monotone functions automatically have bounded variation on closed
intervals, and conversely functions of bounded variation can be
expressed as differencesf of monotone increasing functions.  However,
an example of Cantor shows that a continuous monotone function can
have derivative equal to $0$ almost everywhere without being constant.
This does not happen with the additional condition of absolute
continuity, which is automatically satisfied by Lipschitz functions.

\section{Calculus, 2}
\label{calculus, 2}
\setcounter{equation}{0}

        Let $(M, d(x, y))$ be a metric space, and let $f$ be a
real-valued function on $M$.  There is also a simple version of
$|\nabla f(p)| \le k$ that can be defined on $M$ for each $k \ge 0$,
which is that
\begin{equation}
\label{limsup_{x to p} frac{|f(x) - f(p)|}{d(x, p)} le k}
        \limsup_{x \to p} \frac{|f(x) - f(p)|}{d(x, p)} \le k.
\end{equation}
To be more explicit, this means that for each $\epsilon > 0$ there is
a $\delta > 0$ such that
\begin{equation}
\label{frac{|f(x) - f(p)|}{d(x, p)} < k + epsilon}
        \frac{|f(x) - f(p)|}{d(x, p)} < k + \epsilon
\end{equation}
for every $x \in M$ with $d(x, p) < \delta$.  If $M$ is ${\bf R}^n$
with the standard metric and $f$ is differentiable at $p$, then this
is equivalent to $|\nabla f(p)| \le k$ in the usual sense.

        Consider the case of a continuous function $f$ on a closed
interval $[a, b]$ in the real line that satisfies (\ref{limsup_{x to
p} frac{|f(x) - f(p)|}{d(x, p)} le k}) at every point in $(a, b)$.  If
$f$ is differentiable at every point in $(a, b)$, then the mean value
theorem implies that $f$ is Lipschitz on $[a, b]$ with constant $k$.
This can be extended to deal with functions that may not be
differentiable, as follows.  If $r$ is a real number such that $|r| >
k$, then one can check that $f(x) + r \, x$ cannot have any local
maxima or minima in $(a, b)$.  This is analogous to the fact that the
derivative of a function is equal to $0$ at a local minimum or maximum
when it exists.  Thus $f(x) + r \, x$ attains its maximum and minimum
on $[a, b]$ at the endpoints, and one can use this to show that
\begin{equation}
        |f(b) - f(a)| \le k \, (b - a).
\end{equation}
The same argument works on any closed subinterval of $[a, b]$, which
implies that $f$ is Lipschitz with constant $k$ on $[a, b]$.

        Suppose now that $p : [a, b] \to M$ is a nice curve, and that
$f : M \to {\bf R}$ satisfies (\ref{limsup_{x to p} frac{|f(x) -
f(p)|}{d(x, p)} le k}) at each point in $p([a, b])$.  In this case,
$f(p(x))$ satisfies the conditions described in the previous
paragraph, and hence is Lipschitz with constant $k$ on $[a, b]$.
Alternatively, suppose that $p : [a, b] \to M$ is a continuous path,
and that (\ref{limsup_{x to p} frac{|f(x) - f(p)|}{d(x, p)} le k})
holds with uniform convergence on $p([a, b])$.  The latter means that
for every $\epsilon > 0$ there is a $\delta > 0$ such that
(\ref{frac{|f(x) - f(p)|}{d(x, p)} < k + epsilon}) holds at every
point in $p([a, b])$.  As in the previous section, one can show
directly that $f(p(x))$ has finite length less than or equal to $k$
times the length of $p$, by estimating sums associated to sufficiently
fine partitions of $[a, b]$.

        If $f$ is Lipschitz with constant $k$, then (\ref{limsup_{x to
p} frac{|f(x) - f(p)|}{d(x, p)} le k}) obviously holds at every point
in $M$.  Lipschitz functions on ${\bf R}^n$ are also differentiable
almost everywhere, as mentioned previously.  Conversely, if a
Lipschitz function $f$ on ${\bf R}^n$ satisfies $|\nabla f(p)| \le k$
almost everywhere, then it can be shown that $f$ is Lipschitz with
constant $k$ on ${\bf R}^n$.

\section{Happy fractals}
\label{happy fractals}
\setcounter{equation}{0}

        Let $(M, d(x, y))$ be a metric space, and suppose that there
is a real number $L \ge 1$ such that for each $x, y \in M$ there is
a nice curve $p : [a, b] \to M$ such that $p(a) = x$, $p(b) = y$, and
\begin{equation}
\label{b - a le L d(x, y)}
        b - a \le L \, d(x, y).
\end{equation}
This means that every pair of points $x$, $y$ in $M$ can be connected
by a curve of length less than or equal to $L \, d(x, y)$.  Of course,
any curve in $M$ connecting $x$ and $y$ has length at least $d(x, y)$.
If $M$ is ${\bf R}^n$ with the standard Euclidean metric, then this
condition holds with $L = 1$, because every pair of points can be
connected by a straight line segment.

        Suppose that $f$ is a real-valued function on $M$ and $k$ is a
nonnegative real number such that (\ref{limsup_{x to p} frac{|f(x) -
f(p)|}{d(x, p)} le k}) holds at every point in $M$.  If $p : [a, b]
\to M$ is a nice curve in $M$, then
\begin{equation}
\label{|f(p(a)) - f(p(b))| le k (b - a)}
        |f(p(a)) - f(p(b))| \le k \, (b - a),
\end{equation}
as in the previous section.  Thus the geometric condition about nice
curves implies that $f$ is Lipschitz with constant $k \, L$ on $M$.
In particular, $f$ is constant when $k = 0$.

        It is not too difficult to see that some standard fractals
like Sierpinski gaskets and carpets and Menger sponges have this
property.  The edges of the various triangles, squares, or cubes can
be used as building blocks for nice curves.  One can use sequences of
these edges at different scales to connect any point to a vertex at
some level of the construction.  However, one should be careful about
points that are close to each other but may not be in the same small
triangle, square, or cube.  They may be in adjacent small triangles,
squares, or cubes, so that one should connect from one to the other.

        Let $M$ be a compact connected smooth submanifold of ${\bf
R}^n$, with the restriction of the standard Euclidean metric to $M$ as
the metric on $M$.  In this case, $M$ certainly satisfies the geometric
condition described before for some $L \ge 1$.  The size of $L$ reflects
the relationship between the intrinsic and extrinsic geometry of $M$.

\section{Happy fractals, 2}
\label{happy fractals, 2}
\setcounter{equation}{0}

        In some situations, distances on a metric space might be
defined in terms of minimization of lengths of paths, so that the
geometric condition discussed in the previous section holds
automatically.  On a Riemannian manifold, for instance, the length of
a smooth path is defined directly in terms of the Riemannian metric,
as the integral of the lengths of the tangent vectors along the path.
The distance between two points is then defined to be the infimum of
the lengths of the paths connecting them.

        To compensate for this, one can ask for other geometric
conditions that limit the complexity of the metric space.  As in
\cite{c-w-1, c-w-2}, one such condition is the \emph{doubling
condition}, that every ball of radius $r$ be contained in the union of
a bounded number of balls of radius $r/2$.  It is easy to see that
Euclidean spaces satisfy this condition, and that subsets of doubling
metric spaces are also doubling with respect to the induced metric.
Thus the doubling condition holds automatically for fractal sets in
${\bf R}^n$, but the condition in the previous section does not.

        \emph{Sub-Riemannian spaces} are another very interesting
class of examples.  One can start with a smooth manifold as before,
but now the paths are restricted to have their tangent vectors in
certain subspaces of the tangent bundle.  Under suitable conditions,
it is still possible to connect any two elements of the manifold by
such a path.  The lengths of these paths can be defined by integrals
in the usual way, and the infimum of the lengths of these paths
between two points determines a metric on the manifold.  Under
suitable conditions again, the resulting metric is compatible with the
standard topology on the manifold and doubling.  However, it is also
fractal.  More precisely, the Hausdorff dimension is still an integer,
but it is larger than the topological dimension of the manifold.

        In addition to a doubling condition on the metric space, there
is often an interesting \emph{doubling measure}, which is a positive
Borel measure such that the measure of a ball is bounded by a constant
times the measure of the ball with the same center and one-half the
radius.  For example, Lebesgue measure has this property on ${\bf
R}^n$, and one can show that the existence of a doubling measure on a
metric space implies that the space is doubling.  A lot of analysis on
Euclidean spaces can be extended to metric spaces with doubling
measures, or \emph{spaces of homogeneous type}, as in \cite{c-w-1,
c-w-2}.  This applies to disconnected spaces like Cantor sets as well
as snowflake spaces and spaces with a lot of rectifiable curves.  With
the latter, one can do more, as in the previous sections.

\section{Very happy fractals}
\label{very happy fractals}
\setcounter{equation}{0}

        On the real line, the fundamental theorem of calculus permits
the behavior of a function to be analyzed in terms of integrals of its
derivative, instead of uniform bounds for the derivative.  There are
also versions of this in arbitrary dimensions, including Poincar\'e
and Sobolev inequalities.  There are analogous results for nilpotent
Lie groups and sub-Riemannian spaces as well.  This is related to
isoperimetric inequalities too, in which the size of a region is
estimated in terms of the size of its boundary.

        In \cite{h-s}, Juha Heinonen and I asked whether some type of
connection with approximately Euclidean or sub-Riemannian geometry was
necessary in order to have analytic or geometric properties like
these, under suitable conditions.  Remarkable examples of Bourdon and
Pajot \cite{b-p} and Laakso \cite{tl} show that there are quite
different spaces with similar features.  These examples have
topological dimension equal to $1$, but different geometry from the
classical Sierpinski gaskets and carpets and Menger sponges.

        As discussed in \cite{s5}, the examples of Bourdon, Pajot, and
Laakso show that stronger hypotheses are required to get positive results.
In a sense, the necessary conditions in H\"ormander's theorem on sums
of squares of vector fields \cite{h} take a step in this direction,
except that one is already working on a smooth manifold by hypothesis.

        The answer seems to lie somewhere between \cite{b-v} and
\cite{ch}.  In the former, regularity of the tangent cone of a metric
space leads to local models based on nilpotent Lie groups.  In the
latter, Poincar\'e inequalities are used to obtain a measurable
tangent bundle related to directional derivatives almost everywhere.
On an ordinary smooth manifold, the tangent space at a point is both
the domain for directional derivatives and a local model for the
geometry of the manifold.  However, these are not quite the same
already for standard sub-Riemannian spaces.

\end{document}